\documentclass[preprint,12pt]{elsarticle}
\biboptions{numbers,sort&compress}
\usepackage{amssymb}
\usepackage{mathrsfs}
\usepackage{amsmath}
\allowdisplaybreaks[4]
\usepackage{hyperref}
\hypersetup{colorlinks=true}
\usepackage{amsthm}
\numberwithin{equation}{section}
\newtheorem{theorem}{Theorem}[section]
\newtheorem{definition}[theorem]{Definition}
\newtheorem{lemma}[theorem]{Lemma}
\newtheorem{remark}[theorem]{Remark}

\usepackage{amssymb}

\journal{arXiv}

\begin{document}

\begin{frontmatter}



\title{Variational approach to the existence of solutions for non-instantaneous impulsive differential equations with perturbation}


\author[Yao]{Wangjin Yao}
\address[Yao]{School of Mathematics and Finance,  Putian University,  Putian, 351100, P.R. China}
\author[Dong]{Liping Dong}
\address[Dong]{College of Mathematics and Informatics,  Fujian Normal University,  Fuzhou, 350117, P.R. China}
\author[Zeng]{Jing Zeng\corref{cor}}
\address[Zeng]{College of Mathematics and Informatics,  Fujian Key Laboratory of Mathematical Analysis and Applications (FJKLMAA), Fujian Normal University, Fuzhou, 350117, P.R. China}
\cortext[cor]{Corresponding author,  email address: zengjing@fjnu.edu.cn. The author is supported by the National Science Foundation of China (Grant No. 11501110) and Fujian Natural Science Foundation (Grant No. 2018J01656).}



\begin{abstract}
In this paper, we study the existence of solutions for second-order non-instantaneous impulsive differential equations with a perturbation term. By variational approach, we obtain the problem has at least one solution under assumptions that the nonlinearities are super-quadratic at infinity, and sub-quadratic at the origin.
\end{abstract}


\end{frontmatter}


\section{Introduction}
\label{}

In this paper, we consider the following problem:
\begin{equation}\label{eq1}
\left\{ {\begin{array}{l}
-u''(t)=D_{x}F_{i}(t,u(t)-u(t_{i+1}))+p(t),\quad t\in(s_{i},t_{i+1}],~i=0,1,2,...,N,\\
u'(t)=\alpha_{i},\qquad \qquad \qquad \qquad \qquad \qquad ~~\quad t\in(t_{i},s_{i}],~i=1,2,...,N,\\
u'(s_{i}^{+})=u'(s_{i}^{-}),\qquad \qquad \qquad \qquad \qquad~ \quad i=1,2,...,N,\\
u(0)=u(T)=0, u'(0)=\alpha_{0},
\end{array}} \right.
\end{equation}
where $0=s_{0}<t_{1}<s_{1}<t_{2}<s_{2}<...<t_{N}<s_{N}<t_{N+1}=T$. For the impulses start abruptly at the points $t_{i}$ and keep the derivative constant on a finite time interval $(t_{i},s_{i}]$, we set $u'(s_{i}^{\pm})=\lim_{s\rightarrow s_{i}^{\pm}}u'(s)$. $\alpha_{i} \ (i=1,...,N)$  are constants, $p(t):(s_{i},t_{i+1}]\rightarrow \mathbb{R}$ belongs to $ L^{2}(s_{i},t_{i+1}] (i=1, ..., N)$.

\smallskip

The mathematical model of real world phenomena, in which discontinuous jump occurs, leads to the impulsive differential equations. The non-instantaneous impulsive differential equation is related to the hemodynamical equilibrium. Hence, it is important to study the non-instantaneous impulsive differential equations with a perturbation term, such as $p(t)$ in \eqref{eq1}. %
As far as we know, the introduction of equation \eqref{eq1} was initiated by Hern$\acute{a}$ndez and O'Regan in \cite{8}. 
In \eqref{eq1}, the action starts abruptly at points $t_{i}$, and remains during a finite time interval
$(t_{i}, s_{i}]$. Obviously, it is a natural generalization of the following classical instantaneous impulsive differential equation:
\begin{equation}\label{eq55}
\left\{ {\begin{array}{l}
-u''(t)=f(t,u(t)),\quad t\in([0, T],\\
u'(t_{i}^{+})-u'(t_{i}^{-})=I_i(u(t_i)),\qquad i=1,2,...,N,\\
u(0)=u(T)=0. 
\end{array}} \right.
\end{equation}
Many classical methods can be used to study the non-instantaneous impulsive differential equations, such as theory of Analytic Semigroup, Fixed-Point theory \cite{6,7,12,13} and so on. For some recent works on this type equation, we refer the readers to \cite{1,4,5,10,11,15,16,17}.


\smallskip

To the best of our knowledge, Variational Method can be used to study some impulsive differential equation. 
Bai-Nieto \cite{2} studied the following linear problem, and obtained the existence and uniqueness of weak solutions.
\begin{equation*}\label{eq2}
\left\{ {\begin{array}{l}
-u''(t)=\sigma_{i}(t),\quad t\in(s_{i},t_{i+1}], i=0,1,2,...,N,\\
u'(t)=\alpha_{i},\quad t\in(t_{i},s_{i}], i=1,2,...,N,\\
u'(s_{i}^{+})=u'(s_{i}^{-}),\quad i=1,2,...,N,\\
u(0)=u(T)=0 , u'(0)=\alpha_{0},
\end{array}} \right.
\end{equation*}
where $\sigma_{i}\in L^{2}((s_{i},t_{i+1}),\mathbb{R})$, $\alpha_{i} \ (i=0,...,N)$  are constants. By Variational Method, Bai-Nieto-Wang \cite{3} obtained at least two distinct nontrivial weak solutions of problem:
\begin{equation*}\label{eq3}
\left\{ {\begin{array}{l}
-u''(t)=D_{x}F_{i}(t,u(t)-u(t_{i+1})),\quad t\in(s_{i},t_{i+1}], i=0,1,2,...,N,\\
u'(t)=\alpha_{i},\quad t\in(t_{i},s_{i}], i=1,2,...,N,\\
u'(s_{i}^{+})=u'(s_{i}^{-}),\quad i=1,2,...,N,\\
u(0)=u(T)=0 , u'(0)=\alpha_{0},
\end{array}} \right.
\end{equation*}
where $D_{x}F_{i}(t,x)$ are the derivatives of $F_{i}(t,x)$ with respect to $x$, $i=0,1,2,...,N.$ Zhang-Yuan \cite{18} considered the following equation with a perturbation term $p(t)$, and obtained infinitely many weak solutions.
\begin{equation*}\label{eq4}
\left\{ {\begin{array}{l}
-u''(t)+\lambda u(t)=f(t,u(t))+p(t), \quad a.e.~t\in[0,T],\\
 \bigtriangleup u'(t_{i})=I_{i}(u(t_{i})),\quad i=1,...,N,\\
 u(0)=u(T)=0,
\end{array}} \right.
\end{equation*}
where $f: [0, T]\times\mathbb{R}\rightarrow \mathbb{R}$ is continuous, the impulsive functions $I_{i}:\mathbb{R}\rightarrow \mathbb{R} (i=1, 2, . . . ,N)$ are continuous and $p(t):[0, T]\rightarrow \mathbb{R}$ belongs to $L^{2}[0, T]$.

\medskip

Motivated by the work of \cite{2,3,18}, we 
obtain the weak solution of the problem \eqref{eq1} by Variational Method. Our main result is a natural extension of \cite{3}. 
We denotes $D_{x}F_{i}(t,x)$ the derivatives of $F_{i}(t,x)$ with respect to $x (i=0, 1, ..., N)$. $ F_{i}(t,x) $ is measurable in $t$ for every $x\in \mathbb{R}$ and continuously differentiable in $x$ for $a.e.\ t\in (s_{i},t_{i+1}]$.

We assume that $\lambda_{1}$ is the first eigenvalue of:
\begin{equation}\label{eq5}
\left\{ {\begin{array}{l}
 \displaystyle -u''(t)=\lambda u(t), \quad t\in[0,T],\\
             \displaystyle u(0)=u(T)=0.
\end{array}} \right.
\end{equation}

\smallskip

Our assumptions are:
\begin{description}
  \item[$(H1)$] There exist $\alpha \in C(\mathbb{R}^{+},\mathbb{R}^{+})$  and $ b \in L^{1}(s_{i},t_{i+1};\mathbb{R}^{+})$ such that
$$|F_{i}(t,x)|\leq \alpha(|x|)b(t),~|D_{x}F_{i}(t,x)|\leq \alpha(|x|)b(t),$$ for all $x\in \mathbb{R}$,
 where $F_{i}(t,0)=0$ for $a.e.\ t\in(s_{i},t_{i+1}) ~(i=0,1,2,...,N)$.
  \item[$(H2)$] There exist constants $\mu_{i}>2$ such that $0<\mu_{i}F_{i}(t,x)\leq xD_{x}F_{i}(t,x)$
 for $a.e.\ t\in (s_{i},t_{i+1}], ~x\in \mathbb{R}\backslash \{0\} (i=0,1,2,...,N).$
  \item[$(H3)$] There exist constant $M$ such that $\sum\limits_{i=0}^{N}\sqrt{t_{i+1}-s_{i}}\|p\|_{L^{2}}<M,$

where $M=\frac{1}{8\beta^{2}}-\frac{1}{2}\sum\limits_{i=1}^{N}|\alpha_{i-1}-\alpha_{i}|-\sum\limits _{i=0}^{N}\int_{s_{i}}^{t_{i+1}}M_{i}(t)dt,~\beta=(T\lambda_{1})^{-\frac{1}{2}}+T^{\frac{1}{2}},~M_{i}(t):=\max \limits_{|x|=1}F_{i}(t,x)~(i=0,1,2,...,N).$
\end{description}
\begin{remark}
  $M$ in $(H3)$ is originated from the proof of Theorem \ref{th1}.
\end{remark}
\begin{theorem} \label {th1}
Suppose that $(H1)$-$(H3)$ hold, then problem \eqref{eq1} has at least one weak solution.
\end{theorem}

The article is organized as following: In Section 2, we present some basic knowledge and preliminary results. In Section 3, we prove Theorem \ref{th1}.

\section{Preliminaries}

In this section, we present some preliminary results which will be used in the proof of our result.

\begin{definition}\label{del}{\bf (\cite{9}, (PS) condition)}
Let $E$ be a real Banach space and $I\in C^{1}(E, \mathbb{R})$. $I$ is said to be satisfying the Palais-Smale condition on $E$ if
any sequence $\{u_{k}\}\in E$ for which $I(u_{k})$ is bounded and $I'(u_{k})\rightarrow0$ as $k\rightarrow\infty$ possesses
a convergent subsequence in $E$.
\end{definition}

\begin{theorem}\label{th2}{\bf(\cite{14}, Mountain Pass Theorem)}
 Let $E$ be a real Banach space and $I\in C^{1}(E,\mathbb{R})$ satisfy the $(PS)$ condition with $I(0)=0$. If $I$ satisfies the following conditions:
\begin{description}
  \item[$(1)$] there exist constants $\rho,\alpha >0$, such that $I|_{\partial B_{\rho}}\geq \alpha$;
  \item[$(2)$] there exists an $e\in E\backslash B_{\rho}$, such that $I(e)\leq 0$,
\end{description}
then $I$ possesses a critical value $c\geq \alpha$. Moreover, $c$ is characterized as $$c=\inf \limits_{g\in \Gamma}\max \limits_{s\in [0,1]}I(g(s)),$$
where $$\Gamma=\{g\in C([0,T],E)|~g(0)=0,g(1)=e\}.$$
\end{theorem}

Next, we introduce the well-known Poincar$\acute{e}$ inequality $$\int_{0}^{T}|u|^{2}dt\leq\frac{1}{\lambda_{1}}\int_{0}^{T}|u'|^{2}dt, ~u\in H_{0}^{1}(0,T),$$
where $\lambda_{1}$ is given in \eqref{eq5}.

\smallskip

In the Sobolev space $H_{0}^{1}(0,T)$, we consider the inner product $(u,v)=\int_{0}^{T}u'(t)v'(t)dt,$
which induces the norm $\|u\|=\left(\int_{0}^{T}|u'(t)|^{2}\right)^{\frac{1}{2}}.$ In $L^{2}[0,T]$ and $C[0,T]$, we define the norms:
$$\|u\|_{L^{2}}=\left(\int_{0}^{T}|u(t)|^{2}dt\right)^{\frac{1}{2}},~~\|u\|_{\infty}= \max\limits_{t\in [0,T]}|u(t)|.$$

By the Mean Value Theorem and the H$\ddot{o}$lder inequality, for any $u\in H_{0}^{1}(0,T)$, we have
\begin{equation}\label{eq106}
\|u\|_{\infty}\leq\beta\|u\|,
\end{equation}
where $\beta=(T\lambda_{1})^{-\frac{1}{2}}+T^{\frac{1}{2}},$ $\lambda_{1}$ is given in \eqref{eq5}.

\smallskip

Take $v\in H_{0}^{1}(0,T)$, multiply \eqref{eq1} by $v$ and integrate from $0$ to $T$, we obtain
\begin{equation*}\label{eq6}
  \begin{split}
\int_{0}^{T}u''vdt=&\int_{0}^{t_{1}}u''vdt+\sum\limits_{i=1}^{N}\int_{t_{i}}^{s_{i}}u''vdt+\sum\limits_{i=1}^{N-1}\int_{s_{i}}^{t_{i+1}}u''vdt+\int_{s_{N}}^{T}u''vdt\\
=&-\int_{0}^{T}u'v'dt+\sum\limits_{i=1}^{N}[u'(t_{i}^{-})-u'(t_{i}^{+})]v(t_{i})+\sum\limits_{i=1}^{N}[u'(s_{i}^{-})-u'(s_{i}^{+})]v(s_{i}).
  \end{split}
\end{equation*}

By \eqref{eq1},

\begin{equation}\label{eq7}
  \begin{split}
 \int_{0}^{T}u''vdt=&-\int_{0}^{T}u'v'dt+\sum\limits_{i=1}^{N}[\alpha_{i-1}-\alpha_{i}]v(t_{i})\\
&-\sum\limits_{i=0}^{N-1}\int_{s_{i}}^{t_{i+1}}(D_{x}F_{i}(t,u(t)-u(t_{i+1}))+p(t))dt)v(t_{i+1}).
  \end{split}
\end{equation}
On the other hand,

\begin{equation}\label{eq8}
  \begin{split}
\int_{0}^{T}u''vdt=&-\sum\limits_{i=0}^{N}\int_{s_{i}}^{t_{i+1}}(D_{x}F_{i}(t,u(t)-u(t_{i+1}))+p(t))vdt+\sum\limits_{i=1}^{N}\int_{t_{i}}^{s_{i}}\frac{d}{dt}[\alpha_{i}]vdt\\
=&-\sum\limits_{i=0}^{N}\int_{s_{i}}^{t_{i+1}}(D_{x}F_{i}(t,u(t)-u(t_{i+1}))+p(t))vdt.
  \end{split}
\end{equation}
Thus, it follows $v(t_{N+1})=v(T)=0$, \eqref{eq7} and \eqref{eq8} that

\begin{equation}\label{eq9}
  \begin{split}
-\int_{0}^{T}u'v'dt+\sum\limits_{i=1}^{N}[\alpha_{i-1}-\alpha_{i}]v(t_{i})=&-\sum\limits_{i=0}^{N}\int_{s_{i}}^{t_{i+1}}(D_{x}F_{i}(t,u(t)-u(t_{i+1}))\\
&+p(t))(v(t)-v(t_{i+1})dt.
 \end{split}
\end{equation}

A weak solution to \eqref{eq1} is a function $u\in H_{0}^{1}(0,T)$ such that \eqref{eq9} holds for any $v\in H_{0}^{1}(0,T)$.

Consider the functional $I:~H_{0}^{1}(0,T)\rightarrow \mathbb{R},$
\begin{equation}\label{eq10}
  \begin{split}
I(u)=&\displaystyle\frac{1}{2}\int_{0}^{T}|u'|^{2}dt-\sum\limits_{i=1}^{N}(\alpha_{i-1}-\alpha_{i})u(t_{i})\\
&-\sum\limits_{i=0}^{N}\int_{s_{i}}^{t_{i+1}}p(t)(u(t)-u(t_{i+1}))dt-\sum\limits_{i=0}^{N}\varphi_{i}(u),
 \end{split}
\end{equation}
where $\varphi_{i}(u):=\displaystyle\int_{s_{i}}^{t_{i+1}}F_{i}(t,u(t)-u(t_{i+1}))dt.$

For $u$ and $v$ fixed in $H_{0}^{1}(0,T)$ and $\lambda\in[-1, 1]$, by \eqref{eq106}, we have
\begin{equation}\label{eq50}
|u(t)-u(t_{i+1})|\leq2\|u\|_{\infty}\leq2\beta\|u\|.
\end{equation}
Hence
$$|u(t)-u(t_{i+1})+\lambda\theta(v(t)-v(t_{i+1}))|\leq2\beta(\|u\|+\|v\|),~\text{for} ~\theta\in(0,1),$$
and for $a.e.$ $t\in (s_{i},t_{i+1}]$,
\begin{align*}
  \begin{split}
&\lim\limits_{\lambda\rightarrow0}\frac{1}{\lambda}\left[F_{i}(t,u(t)-u(t_{i+1})+\lambda(v(t)-v(t_{i+1})))-F_{i}(t,u(t)-u(t_{i+1}))\right]\\
=&D_{x}F_{i}(t,u(t)-u(t_{i+1}))(v(t)-v(t_{i+1})).
  \end{split}
\end{align*}

By $(H1)$, \eqref{eq50} and the Mean Value Theorem, we obtain
\begin{equation*}
  \begin{split}
&\left|\frac{1}{\lambda}\left[F_{i}(t,u(t)-u(t_{i+1})+\lambda(v(t)-v(t_{i+1})))-F_{i}(t,u(t)-u(t_{i+1}))\right]\right|\\
=&\bigg|D_{x}F_{i}(t,u(t)-u(t_{i+1})+\lambda\theta(v(t)-v(t_{i+1}))(v(t)-v(t_{i+1}))\bigg|\\
\leq&\max \limits_{z\in[0,2\beta(\|u\|+\|v\|)]}a(z)2\beta\|v\|b(t)\in L^{1}(s_{i},t_{i+1};\mathbb{R}^{+}).
\end{split}
\end{equation*}
Lebesgue's Dominated Convergence Theorem shows that
\begin{equation*}\label{eq11}
(\varphi_{i}'(u),v)=\int_{s_{i}}^{t_{i+1}}D_{x}F_{i}(t,u(t)-u(t_{i+1}))(v(t)-v(t_{i+1}))dt.
\end{equation*}
Moreover, $\varphi_{i}'(u)$ is continuous. So $I\in C^{1}(H_{0}^{1}(0,T),\mathbb{R})$ and
\begin{equation}\label{eq12}
  \begin{split}
I'(u)v=&\int_{0}^{T}u'v'dt+\sum\limits_{i=1}^{N}[\alpha_{i-1}-\alpha_{i}]v(t_{i})\\
&-\sum\limits_{i=0}^{N}\int_{s_{i}}^{t_{i+1}}\left(D_{x}F_{i}(t,u(t)-u(t_{i+1}))+p(t)\right)(v(t)-v(t_{i+1}))dt.
\end{split}
\end{equation}
Then the correspond critical points of $I$ are the weak solutions of the problem \eqref{eq1}.

\begin{lemma} \label{le2} {\bf(\cite{3})}
If assumption $(H2)$ holds, then for each $i=0,1,2,..,N$, there exist $M_{i},m_{i},b_{i}\in L^{1}(s_{i}, t_{i+1})$ which are
almost everywhere positive such that
$$F_{i}(t,x)\leq M_{i}(t)|x|^{\mu_{i}},~for ~a.e.~t\in(s_{i}, t_{i+1}],~and~|x|\leq1,$$
and
$$F_{i}(t,x)\geq m_{i}(t)|x|^{\mu_{i}}-b_{i}(t),~for ~a.e.~t\in(s_{i}, t_{i+1}],~and~x\in\mathbb{R},$$
where $m_{i}(t):=\min \limits_{|x|=1}F_{i}(t,x)$, $M_{i}(t):=\max \limits_{|x|=1}F_{i}(t,x),~a.e.~t\in(s_{i},t_{i+1}].$
\end{lemma}
\begin{remark}
Lemma \ref{le2} implies that $D_{x}F_{i}(t,x)\ (i=1, ..., N)$ are super-quadratic at infinity, and sub-quadratic at the origin.
\end{remark}
\begin{lemma}\label{le3}
Suppose that $(H1)$, $(H2)$ hold, then $I$ satisfies the (PS) condition.
\end{lemma}

\noindent{\bf Proof:} Let $\{u_{k}\}\subset H_{0}^{1}(0,T)$ such that $\{I (u_{k})\}$ be a bounded sequence and $\lim \limits_{k\rightarrow \infty}I'(u_{k})=0$.

\smallskip

By \eqref{eq106},
\begin{equation}\label{eq17}
|\sum\limits_{i=1}^{N}(\alpha_{i-1}-\alpha_{i})u(t_{i})|\leq\sum\limits_{i=1}^{N}|\alpha_{i-1}-\alpha_{i}|\|u\|_{\infty}\leq\sum\limits_{i=1}^{N}|\alpha_{i-1}-\alpha_{i}|\beta\|u\|. \end{equation}
There exists constant $C_{1}>0$ such that $$|I(u_{k})|\leq C_{1},~|I'(u_{k})|\leq C_{1}.$$

First, we prove that $\{u_{k}\}$ is bounded. Let $\mu:=\min\{\mu_{i}:i=0,1,2,...,N\}$, by \eqref{eq10}, \eqref{eq17} and $(H2)$, we obtain
\begin{equation*}
\begin{split}
\int_{0}^{T}|u_{k}'|^{2}dt=& 2I(u_{k})+2\sum\limits_{i=1}^{N}(\alpha_{i-1}-\alpha_{i})u_{k}(t_{i})\\
&+2\sum\limits_{i=0}^{N}\int_{s_{i}}^{t_{i+1}}p(t)(u_{k}(t)-u_{k}(t_{i+1}))dt\\
&+2\sum\limits_{i=0}^{N}\int_{s_{i}}^{t_{i+1}}F_{i}(t,u_{k}(t)-u_{k}(t_{i+1}))dt,\\
\leq& 2C_{1}+2\beta\sum\limits_{i=1}^{N}|\alpha_{i-1}-\alpha_{i}|\|u_{k}\|\\
&+2\sum\limits_{i=0}^{N}\int_{s_{i}}^{t_{i+1}}p(t)(u_{k}(t)-u_{k}(t_{i+1}))dt\\
&+\frac{2}{\mu}\sum\limits_{i=0}^{N}\int_{s_{i}}^{t_{i+1}}D_{x}F_{i}(t,u_{k}(t)-u_{k}(t_{i+1}))(u_{k}(t)-u_{k}(t_{i+1}))dt,
\end{split}
\end{equation*}
which combining \eqref{eq12} yields that
\begin{align*}
(1-\frac{2}{\mu})\|u_{k}\|^{2}\leq&2C_{1}+(2+\frac{2}{\mu})\sum\limits_{i=1}^{N}|\alpha_{i-1}-\alpha_{i}|\beta\|u_{k}\|-\frac{2}{\mu}I'(u_{k})u_{k}\\
&+(2-\frac{2}{\mu})\sum\limits_{i=0}^{N}\int_{s_{i}}^{t_{i+1}}p(t)(u_{k}(t)-u_{k}(t_{i+1}))dt,\\
\leq& 2C_{1}+(2+\frac{2}{\mu})\beta\sum\limits_{i=1}^{N}|\alpha_{i-1}-\alpha_{i}|\|u_{k}\|+\frac{2}{\mu}C_{1}\beta\|u_{k}\|\\
&+2(2-\frac{2}{\mu})\beta\sum\limits_{i=0}^{N}\|u_{k}\|\sqrt{t_{i+1}-s_{i}}\|p\|_{L^{2}}.
\end{align*}

Since $\mu>2$, it follow that $\{u_{k}\}$ is bounded in $H_{0}^{1}(0,T)$.

\smallskip

Therefore, there exists a subsequence also denoted by $\{u_{k}\}\in H_{0}^{1}(0,T)$ such that
\begin{equation*}
\begin{split}
&u_{k}\rightharpoonup u, ~~ \text{in} ~H_{0}^{1}(0,T),\\
&u_{k} \rightarrow u, ~~\text{in} ~L^{2}(0,T),\\
&u_{k} \rightarrow u, ~~\text{uniformly in} ~[0,T],~~\text{as} ~k\rightarrow \infty.
\end{split}
\end{equation*}

Since
\begin{equation*}
\begin{split}
	|u_{k}(t)-u_{k}(t_{i+1})-u(t)+u(t_{i+1})|\leq& |u_{k}(t)-u(t)|+|u(t_{i+1})-u_{k}(t_{i+1})|\\
	\leq& 2\|u_{k}-u\|\rightarrow 0, \quad \text{as}~~ k\rightarrow \infty.
\end{split}
\end{equation*}

Hence
\begin{equation*}
\begin{split}
\sum\limits_{i=0}^{N}\int_{s_{i}}^{t_{i+1}}&(D_{x}F_{i}(t,u_{k}(t)-u_{k}(t_{i+1}))-D_{x}F_{i}(t,u(t)-u(t_{i+1})))\\
\cdot&(u_{k}(t)-u_{k}(t_{i+1})-u(t)+u(t_{i+1}))dt\rightarrow 0,
\end{split}
\end{equation*}

$$|\langle I'(u_{k})-I'(u),u_{k}-u\rangle|\leq\|I'(u_{k})-I'(u)\|\|u_{k}-u\|\rightarrow 0.$$

Moreover, we obtain
\begin{equation*}
\begin{split}
&\langle I'(u_{k})-I'(u),u_{k}-u\rangle\\
=&\|u_{k}-u\|-\sum\limits_{i=0}^{N}\int_{s_{i}}^{t_{i+1}}(D_{x}F_{i}(t,u_{k}(t)-u_{k}(t_{i+1}))\\
&-D_{x}F_{i}(t,u(t)-u(t_{i+1})))(u_{k}(t)-u_{k}(t_{i+1})-u(t)+u(t_{i+1}))dt,
\end{split}
\end{equation*}
so $\|u_{k}-u\|\rightarrow 0$ as $k\rightarrow +\infty$. That is, $\{u_{k}\}$ converges strongly to $u$ in $H_{0}^{1}(0,T)$.

\smallskip

Thus, $I$ satisfies the (PS) condition. $\Box$

\section{Proof of theorem}

\noindent{\bf Proof of Theorem \ref{th1}} We found that $I(0)=0$ and $I\in C^{1}(H_{0}^{1}(0,T),\mathbb{R})$. By Lemma \ref{le3}, we obtain $I$ satisfies (PS) condition. By Lemma
 \ref{le2} and \eqref{eq50}, we have
 \begin{equation*}
 \begin{split}
\int_{s_{i}}^{t_{i+1}}F_{i}(t,u(t)-u(t_{i+1})dt&\leq\int_{s_{i}}^{t_{i+1}}M_{i}(t)|u(t)-u(t_{i+1})|^{\mu_{i}}dt\\
 &\leq\int_{s_{i}}^{t_{i+1}}M_{i}(t)|2\beta\|u\||^{\mu_{i}}dt,
 \end{split}
 \end{equation*}
and
$$\sum\limits_{i=1}^{N}(\alpha_{i-1}-\alpha_{i})u(t_{i})\leq\sum\limits_{i=1}^{N}|\alpha_{i-1}-\alpha_{i}|\|u\|_{\infty}\leq\sum\limits_{i=1}^{N}|\alpha_{i-1}-\alpha_{i}|\beta\|u\|,$$
$$\sum\limits_{i=0}^{N}\int_{s_{i}}^{t_{i+1}}p(t)(u(t)-u(t_{i+1}))\leq\sum\limits_{i=0}^{N}2\beta\|u\|\sqrt{t_{i+1}-s_{i}}\|p\|_{L^{2}}.$$
By \eqref{eq10},
\begin{equation}\label{eq19}
\begin{split}
I(u)\geq&\frac{1}{2}\|u\|^{2}-\sum\limits_{i=1}^{N}|\alpha_{i-1}-\alpha_{i}|\beta\|u\|-\sum\limits_{i=0}^{N}\int_{s_{i}}^{t_{i+1}}M_{i}(t)|2\beta\|u\||^{\mu_{i}}dt\\
&-\sum\limits_{i=0}^{N}2\beta\|u\|\sqrt{t_{i+1}-s_{i}}\|p\|_{L^{2}}.
\end{split}
\end{equation}
Take $\|u\|=\frac{1}{2\beta}$, then $|u(t)-u(t_{i+1})|\leq 1$, so
\begin{equation*}
  \begin{split}
&\sum\limits_{i=1}^{N}|\alpha_{i-1}-\alpha_{i}|\beta\|u\|\leq\frac{1}{2}\sum\limits_{i=1}^{N}|\alpha_{i-1}-\alpha_{i}|,\\
&\sum\limits_{i=0}^{N}\int_{s_{i}}^{t_{i+1}}M_{i}(t)|2\beta\|u\||^{\mu_{i}}dt\leq\sum\limits _{i=0}^{N}\int_{s_{i}}^{t_{i+1}}M_{i}(t)dt,\\
&\sum\limits_{i=0}^{N}2\beta\|u\|\sqrt{t_{i+1}-s_{i}}\|p\|_{L^{2}}\leq\sum\limits_{i=0}^{N}\sqrt{t_{i+1}-s_{i}}\|p\|_{L^{2}}.
   \end{split}
\end{equation*}
Hence,
\begin{align*}
I(u)=&\frac{1}{2}\|u\|^{2}-\sum\limits_{i=1}^{N}(\alpha_{i-1}-\alpha_{i})u(t_{i})\\
&-\sum\limits_{i=0}^{N}\int_{s_{i}}^{t_{i+1}}p(t)(u(t)-u(t_{i+1}))dt-\sum\limits_{i=0}^{N}\int_{s_{i}}^{t_{i+1}}F_{i}(t,u(t)-u(t_{i+1}))dt,\\
\geq&\frac{1}{2}\|u\|^{2}-\frac{1}{2}\sum\limits_{i=1}^{N}|\alpha_{i-1}-\alpha_{i}|-\sum\limits _{i=0}^{N}\int_{s_{i}}^{t_{i+1}}M_{i}(t)dt-\sum\limits_{i=0}^{N}\sqrt{t_{i+1}-s_{i}}\|p\|_{L^{2}}.
\end{align*}

By $(H3)$, $I(\frac{1}{2\beta})>0$ and satisfies the condition (1) in Theorem \ref{th2}. Let $\xi>0$ and $w\in H_{0}^{1}(0,T)$ with $\|w\|=1$. We can see that $w(t)$ is not a constant for $a.e.~[s_{i}, t_{i+1}]$. By Lemma \ref{le2},
\begin{equation*}
\begin{split}
\int_{s_{i}}^{t_{i+1}}F_{i}(t,(w(t)-w(t_{i+1}))\xi)dt\geq & \left(\int_{s_{i}}^{t_{i+1}}m_{i}(t)|w(t)-w(t_{i+1})|^{\mu_{i}}dt\right)\xi^{\mu_{i}}\\&-\int_{s_{i}}^{t_{i+1}}b_{i}(t)dt.
\end{split}
\end{equation*}

Let $W_{i}:=\int_{s_{i}}^{t_{i+1}}m_{i}(t)|w(t)-w(t_{i+1})|^{\mu_{i}}dt$, then
 $$0\leq W_{i}\leq(2\beta)^{\mu_{i}}\int_{s_{i}}^{t_{i+1}}m_{i}(t)dt,~W_{0}\geq0.$$

We can select the interval $[0, t_{1}]$ and prove $w(t)$ is not a constant for $a.e. ~[0, t_{1}]$. In fact, we suppose that $\int_{0}^{t_{1}}m_{0}(t)|w(t)-w(t_{1})|^{\mu_{0}}dt=0$. Since $m_{0}(t)$ is positive, then $w(t)=w(t_{1})$ for $a.e. ~[0,t_{1}]$. A
contradiction with the assumption on $w$.

By \eqref{eq10}, we obtain
\begin{align*}
I(\xi w)=&\frac{1}{2}\xi^{2}w^{2}-\sum\limits_{i=1}^{N}(\alpha_{i-1}-\alpha_{i})w(t_{i})\xi-\sum\limits_{i=0}^{N}\int_{s_{i}}^{t_{i+1}}p(t)(w(t)-w(t_{i+1})\xi)dt\\
&-\sum\limits_{i=0}^{N}\int_{s_{i}}^{t_{i+1}}F_{i}(t,(w(t)-w(t_{i+1}))\xi)dt,\\
\leq&\frac{1}{2}\xi^{2}+\sum\limits_{i=1}^{N}|\alpha_{i-1}-\alpha_{i}|\beta\xi+2\beta\xi\sum\limits_{i=0}^{N}\sqrt{t_{i+1}-s_{i}}\|p\|_{L^{2}}-\sum\limits_{i=0}^{N}W_{i}\xi^{\mu_{i}}\\
&+\sum\limits_{i=0}^{N}\int_{s_{i}}^{t_{i+1}}b_{i}(t)dt.
\end{align*}

Since $\mu_{i}>2$, the above inequation implies that $I(\xi w)\rightarrow -\infty$ as  $\xi \rightarrow \infty$, that is,  there exists a $\xi\in \mathbb{R}\backslash \{0\}$ such that $\|\xi w\|>\frac{1}{2\beta}$ and $I(\xi w)\leq0$. The proof of Theorem \ref{th1} is completed. $\Box$





\newpage
\noindent \bfseries Acknowledgments \mdseries

\smallskip

Jing Zeng is supported by the National Science Foundation of China (Grant No. 11501110) and Fujian Natural Science Foundation (Grant No. 2018J01656).
\bigskip

\medskip
\noindent \bfseries References \mdseries

\end{document}